\documentclass[11pt]{amsart}

\usepackage[T1]{fontenc}
\usepackage{lmodern}
\usepackage{amsmath,amssymb,amsthm}
\usepackage{mathtools}
\usepackage{xcolor}
\usepackage[all,cmtip]{xy}
\usepackage[hidelinks]{hyperref}

\newcommand{\DDD}{\mathcal{D}}

\newcommand{\FFF}{\mathcal{F}}
\newcommand{\GGG}{\mathcal{G}}

\newcommand{\III}{\mathcal{I}}
\newcommand{\JJJ}{\mathcal{J}}

\newcommand{\MMM}{\mathcal{M}}
\newcommand{\NNN}{\mathcal{N}}
\newcommand{\OOO}{\mathcal{O}}

\let\originalAngstrom\AA
\renewcommand{\AA}{\ifmmode\mathbf{A}\else\originalAngstrom\fi}
\let\originalCapitalSharpS\SS
\renewcommand{\SS}{\ifmmode\mathbf{S}\else\originalCapitalSharpS\fi}

\newcommand{\LL}{\mathbf{L}}

\newcommand{\PP}{\mathbf{P}}

\newcommand{\RR}{\mathbf{R}}

\newcommand{\ZZ}{\mathbf{Z}}

\numberwithin{equation}{section}
\newtheorem{theorem}[equation]{Theorem}
\newtheorem{proposition}[equation]{Proposition}
\newtheorem{lemma}[equation]{Lemma}

\theoremstyle{definition}

\theoremstyle{remark}
\newtheorem{remark}[equation]{Remark}
\newtheorem{example}[equation]{Example}

\DeclareMathOperator{\Bl}{Bl}

\DeclareMathOperator{\Tor}{Tor}
\DeclareMathOperator{\Spec}{Spec}
\DeclareMathOperator{\Supp}{Supp}
\DeclareMathOperator{\Sym}{Sym}
\DeclareMathOperator{\gr}{gr}
\DeclareMathOperator{\rank}{rank}

\DeclareMathOperator{\id}{id}

\DeclareMathOperator{\pr}{pr}
\DeclareMathOperator{\Span}{span}
\newcommand{\m}{\mathfrak m}

\newcommand{\chiR}{\chi^R}

\title[Positivity of intersection numbers]
      {A classical proof of positivity of Serre's intersection numbers}
\date{\today}
\author{Eike Lau}
\address{Fakult\"{a}t f\"{u}r Mathematik,
Universit\"{a}t Bielefeld, D-33501 Bielefeld}

\begin{document}
\begin{abstract}
Recently OpenAI proved positivity of Serre's intersection multiplicity. We give a different proof, starting from Gabber's proof of nonnegativity and vanishing via de Jong alterations. As a first step, the multiplicity of cyclic modules is positive if the addition of the associated tangent cones has dense image in the ambient tangent space. In complementary dimensions a symmetric group action on a finite free cover gives this density.
The main argument is AI-generated.
\end{abstract}

\maketitle

\section{Introduction}\label{sec:introduction}

Let $R$ be a regular local ring of dimension $d$ and 
$M,N$ finite $R$-modules such that $M\otimes_R N$ has finite length. 
Serre \cite{Serre} defined the intersection
multiplicity
\[
  \chiR(M,N)=\sum_{i=0}^d(-1)^i
       \ell_R\bigl(\Tor_i^R(M,N)\bigr),
\]
proved the dimension inequality $\dim M+\dim N\le d$, and conjectured that
\[
\begin{aligned}
  \dim M+\dim N=d\ &\Longrightarrow\ \chiR(M,N)>0, 
  \qquad \text{(positivity)}\\
  \dim M+\dim N<d\ &\Longrightarrow\ \chiR(M,N)=0,
  \qquad \text{(vanishing)}
\end{aligned}
\]
which implies $\chiR(M,N)\ge0$ (nonnegativity). 

Serre proved these assertions if $R$ is of equal characteristic or unramified (more generally if the completion of $R$ is a power series ring over a field or over a discrete valuation ring), by reduction to the diagonal and Hilbert--Samuel multiplicity. General vanishing was proved independently by
Roberts \cite{RobertsVanishing,RobertsChern1987} via localized Chern characters and by Gillet--Soul\'e \cite{GilletSoule,GilletSoule1987} via Adams operations in algebraic $K$-theory.
Gabber proved nonnegativity (and vanishing again) using de Jong's regular alterations \cite{deJongAlterations}, as explained in \cite{Berthelot,HochsterGabber}.
We refer to \cite{Hochster1978,RobertsICM1991,Roberts2012,BhattHochsterMa2026} for further results, background, and history.

In October 2026, OpenAI \cite{OpenAI2026} presented a proof of positivity using
perfectoid methods \cite{BhattScholze2022,CaiEtAl2026} and
rational $K$-theoretic descent with support based on \cite{LandTamme2019}.
Ma \cite{Ma2026} simplified the proof by replacing its use of
asymptotic complexes with a reduction to the diagonal via
$A_{\mathrm{inf}}$.

\subsection*{Positivity via alterations}

We give a different proof of positivity, building on Gabber's method. By approximation one can assume that $R$ is \emph{geometric}, i.e.\ essentially of finite type over a complete discrete valuation ring.
Let $k=R/\m$ be the residue field.
We say that $R$ \emph{has a first-order coefficient field} if
the residue map $R/\m^2\to k$ has a ring section.
For an ideal $I\subseteq R$, let
\[
C_{R/I}=\Spec\gr_{\m}(R/I)
\] 
denote its tangent cone in $T_R=\Spec\Sym(\m/\m^2)$.
Analysing tangent directions in the alteration construction yields the following criterion.

\begingroup
\renewcommand{\theequation}{\getrefnumber{prop:cone-criterion-proof}}
\renewcommand{\theHequation}{intro-dominant-cone}
\begin{proposition}
\label{thm:dominant-cone-intro}
Assume that $R$ is geometric with first-order coefficient field and algebraically closed residue field. Let
$I,J\subseteq R$ be ideals such that $R/(I+J)$ has finite length.  
If the addition map
\[
 C_{R/I}\times_k C_{R/J}\longrightarrow T_R,\qquad
 (x,y)\longmapsto x+y
\]
is dominant, then $\chiR(R/I,R/J)>0$.
\end{proposition}
\addtocounter{equation}{-1}
\endgroup

\begin{proof}[Sketch of proof]
Assume that $I=P$ and $J=Q$ are prime and not maximal.
Choose a regular alteration $Z'\to \Spec(R/Q)$ that factors through the blowup at the closed point, and embed $Z'$ into $\PP^N_R$.
By Gabber's criterion, recalled in Proposition~\ref{prop:gabber-numerical} below, $\chiR(R/P,R/Q)$ is positive iff a certain morphism
\[
q:\widetilde C\to T_R\oplus T_{\mathrm v}
\] 
is dominant. Here $\tilde C$ arises from the normal cone of 
$Z'\cap\PP^N_{R/P}$ in $\PP^N_{R/P}$, and $T_{\mathrm v}$ is the space of global vector fields on $\PP^N_k$. 
The following diagram (Proposition~\ref{prop:criterion-incidence-direct}) relates $q$ to the cone addition.
\[
\xymatrix@C+1em@M+.2em{
 L\times C_{R/P}\times T_{\mathrm v}\ar[r] 
    \ar[d] 
 & \widetilde C\ar[d]^q\\
 C_{R/Q}\times C_{R/P}\times T_{\mathrm v}
    \ar[r]^-{+}
 & T_R\oplus T_{\mathrm v}
}
\]
Here $L$ is the normal line bundle of the special fiber $Z'_s\subseteq Z'$, and  the left vertical arrow is surjective. Hence dominant cone addition implies that $q$ is dominant.
\end{proof}

Proposition~\ref{prop:cone-criterion-proof} is weaker than general positivity because complementary dimensions do not give dominant cone addition, but the latter can be arranged by passing to a finite free regular cover of $R$ carrying a symmetric-group action, provided $\dim(R)+1$ is invertible in $R$, which is also easy to achieve. The details appear in \S\ref{sec:symmetric-covers-direct}.

\subsection*{AI disclosure}
A first version of the present proof was generated by ChatGPT Astra in response to a request to simplify the proof in \cite{OpenAI2026} and was substantially revised in subsequent discussions. The author wrote the final manuscript with AI assistance.

\section{Preliminaries on normal cones}
\label{sec:cone-preliminaries-proposal}

For a closed immersion $Z\to X$ with ideal sheaf $\III$,
its normal cone is denoted by
\[
 C_{Z/X}=\Spec_Z\gr_{\III}\OOO_X.
\]
If the immersion is regular, this is the normal bundle $N_{Z/X}$.
For a Noetherian local ring $A$ with maximal ideal $\m$, 
its tangent cone is $C_A=\Spec\gr_{\m}A$.

\begin{remark} 
If $X$ is noetherian, then
\begin{equation}\label{eq:cone-dimension-preliminaries}
 \dim C_{Z/X}=\sup_{z\in Z}\dim\OOO_{X,z}
\end{equation}
by \cite[Theorem~4.5.6(b)]{BrunsHerzog}.
In particular, $\dim C_A=\dim A$.
\end{remark}

\begin{lemma} 
\label{lem:cone-comparison-preliminaries}
Let $S\to Y\to X$ be closed immersions of
Noetherian schemes, and let $f:X'\to X$ be flat of finite type.
Put $S'=S\times_X X'$ and $Y'=Y\times_X X'$.
Let $D\to S'$ be a regular closed immersion such that
$D\to S$ is flat. Then the natural square
\[
\xymatrix@M+0.2em{
 C_{D/Y'}\ar[r]\ar[d] & C_{D/X'}\ar[d]\\
 C_{S/Y}\ar[r] & C_{S/X}
}
\]
is Cartesian.
\end{lemma}

\begin{proof}
Locally along $D$, choose a regular sequence $g_1,\ldots,g_c$
defining $D$ in $S'$ and lift it to functions on $X'$. Put
$F=(f,g_1,\ldots,g_c):X'\to X\times\AA^c$.

We claim that $F$ is flat along $D$.
By the fiberwise flatness criterion \cite[Tag~039C]{StacksProject}, 
this holds iff the restriction $F_S:S'\to S\times\AA^c$ is flat along $D$, which holds by the local flatness criterion  \cite[Tag~00ML]{StacksProject} because $D$ is flat over $S$ and the relevant $\Tor_1$ vanishes due to the regular sequence, using the Koszul resolution.

Since the flat locus is open we may assume that $F$ is flat.
Then flat base change for normal cones \cite[Tag~0635]{StacksProject} reduces the lemma
to the case $X'=X\times\AA^c$ where  $S=D\to S'=S\times\AA^c$ is the zero section. In that case we have $C_{D/X'}=C_{D/X}\times\AA^c$ and $C_{D/Y'}=C_{D/Y}\times\AA^c$ compatibly, and the assertion follows.
\end{proof}

\begin{lemma} 
\label{lem:cone-minimal-primes-preliminaries}
Let $A$ be a noetherian local ring. Then $C_A$ is the union of $C_{A/P_i}$ for the minimal primes $P_i$ of $A$ as a topological space.
\end{lemma}

\begin{proof}
Let $X$ and $X_i$ be the blowup of $A$ and $A/P_i$ at the maximal ideals, $E\subseteq X$ and $E_i\subseteq X_i$ the exceptional fibers. The morphism $\coprod_i X_i\to X$ is surjective since it is proper and the image contains a dense open set. Hence $\coprod_i E_i\to E$ is surjective, and the assertion follows easily by passing from exceptional fibers to cones.
\end{proof}

\section{Approximation}
\label{app:geometric-approximation}

The situation admits a geometric approximation as follows. 

\begin{lemma} 
\label{lem:geometric-approximation-general}
Let $M,N$ be finite modules over a regular local ring $R$ such that  $M\otimes_RN$ has finite length.
There exist a geometric regular local ring $R_0$ with algebraically closed residue field and finite $R_0$-modules $M_0,N_0$ with
$\dim R_0=\dim R$,
$\dim M_0=\dim M$,
$\dim N_0=\dim N$ such that
$M_0\otimes_{R_0}N_0$ has finite length,
and
\[
\chi^{R_0}(M_0,N_0)=\chiR(M,N).
\]
\end{lemma}

This is proved in \cite[\S6.1.2]{Berthelot} and \cite[\S1 and Appendix]{HochsterGabber} assuming that $R$ is ramified of mixed characteristic, the remaining case of Serre's conjectures, but this hypothesis is not essential.\footnote{Lemma \ref{lem:geometric-approximation-general} is also a formal consequence of Serre's conjectures, hence true for non-ramified $R$ by \cite{Serre} without actual approximation.}

\section{Gabber's method}
\label{sec:gabber-method}

We recall Gabber's results following \cite{Berthelot} and \cite{HochsterGabber}, adapted for later use in \S\ref{sec:cone-criterion-proof}.

\subsection{Euler characteristics with proper support}

Let $(R,\m,k)$ be a Noetherian local ring.
For a bounded complex $K$ of $R$-modules
with finite-length cohomology we write
\[
 \chi_R(K)=\sum_i(-1)^i\ell_RH^i(K).
\]
If $\MMM,\NNN$ are coherent sheaves on an $R$-scheme $X$ of finite type such that $\MMM\otimes_{\OOO_X}\NNN$ has proper support lying over the closed point of $\Spec R$ and one of $\MMM,\NNN$
is perfect (locally of finite Tor-amplitude), we write
\[
 \chi_X(\MMM,\NNN)
 =\chi_R\mathbf R\Gamma
       (X,\MMM\otimes_{\OOO_X}^{\mathbf L}\NNN).
\]

\begin{lemma} 
\label{lem:normal-cone-euler-compare}
Let $X$ be separated and of finite type over $R$, let $i:Z\to X$ be a regular
closed immersion with ideal $\III$, and let $\MMM$ be a coherent
sheaf on $X$ such that
$Z\cap\Supp\MMM$ is proper over $R$
and lies over the closed point. 
View $\gr_{\III}\MMM$ as a coherent sheaf on
$E=N_{Z/X}$.
Then
\begin{equation}
\label{eq:normal-cone-euler}
\chi_X(\MMM,\OOO_Z)
 =\chi_E(\gr_{\III}\MMM,\OOO_Z).
\end{equation}
\end{lemma}

\begin{proof}
As explained in \cite[\S6.1.4, p.~299]{Berthelot} and \cite[\S3]{HochsterGabber}, this follows from a global version of Serre's filtered Koszul spectral sequence \cite[Ch.~IV, A.3]{Serre}
where the superscript $p$ denotes the homogeneous degree,
\[
E_1^{p,q}=
\Tor_{-p-q}^{\OOO_E}
 \bigl(\gr_{\III}\MMM,\OOO_Z\bigr)^p
\ \Longrightarrow\
\Tor_{-p-q}^{\OOO_X}(\MMM,\OOO_Z).
\qedhere
\]
\end{proof}

\begin{remark}
For an alternative proof of Lemma \ref{lem:normal-cone-euler-compare}, let
$\DDD_{Z/X}=\Spec_X(\OOO_X[\tau,\III/\tau])$
be the affine deformation to the normal cone.
The regular immersion $i$ induces a regular immersion
$j:Z\times\AA^1\to\DDD_{Z/X}$, whose fibers at $\tau=1$ and
$\tau=0$ are $i:Z\to X$ and the zero section
$Z\to E$ \cite[\S5.1]{FultonIntersection}.
 The extended Rees sheaf
\[
 \widetilde\MMM=
 \bigoplus_{n\in\mathbb Z}\III^{\max(n,0)}\MMM\,\tau^{-n}
\]
has derived fibers $\MMM$ and $\gr_{\III}\MMM$ at these two points.
Consequently, the derived fibers of
\[
 K=\RR\Gamma(Z\times\AA^1,\LL j^*\widetilde\MMM)
\]
compute the two sides of \eqref{eq:normal-cone-euler}.
The cohomology of $K$ is bounded and finite
over $R[\tau]$ and is annihilated by a power of $\m$, and 
\eqref{eq:normal-cone-euler} follows because
for a finite $k[\tau]$-module $L$, the 
complex $L\xrightarrow{\tau-c}L$ has Euler
characteristic $\rank_{k[\tau]}L$ for
every $c\in k$.
\end{remark}

\begin{lemma} 
\label{lem:gabber-generated}
Let $F$ be a vector bundle of rank $r$ over a proper $k$-scheme $D$ and let $v:D\times_k T\to F$ be a surjective vector bundle homomorphism for a finite-dimensional $k$-vector space $T$. 
If a coherent sheaf $\MMM$ on $F$
has support dimension at most $r$, then $(v^*\MMM)_{k(T)}$  over $D_{k(T)}$ has finite support, and
\begin{equation}
\label{Eq:chiF-dimH0}
 \chi_F(\MMM,\OOO_D)
 =\dim_{k(T)} H^0\bigl(D_{k(T)},(v^*\MMM)_{k(T)}\bigr).
\end{equation}
Here $k(T)$ is the function field of $T$.
The value is positive iff $\Supp(v^*\MMM)\to T$ is dominant. In particular, the value is zero if $\dim\Supp\MMM<r$.
\end{lemma}

\begin{proof}
Let $i_0:D\to F$ be the zero section.
Then 
\[
 \chi_{F}(\MMM,\OOO_D)
 =\chi_k\mathbf R\Gamma
       (D,\mathbf Li_0^*\MMM).
\]
The right hand side and the dimension hypothesis are preserved when $F$ and $\MMM$ are replaced by $D\times_kT$ and $v^*\MMM$. Thus assume that $F=D\times_kT$ and let $q:F\to T$ be the projection. The generic fiber of $\Supp\MMM\to T$ is finite by the dimension hypothesis. Moreover $\RR q_*\MMM$ is represented by a perfect complex $K$ over $k[T]$ with generic fiber 
$H^0(D_{k(T)},\MMM_{k(T)})$ and zero fiber $\RR\Gamma(D,\LL i_0^*\MMM)$. 
Hence both sides of \eqref{Eq:chiF-dimH0} coincide with the alternating rank of $K$.
\end{proof}

\subsection{The alteration and its normal cone}
\label{subsec:alteration-normal-cone}

Now let $(R,\m,k)$ be a regular local ring,
and assume that $R$ is geometric with algebraically closed
residue field and $R$ has a first-order coefficient field;
choose a section $k\to R/\m^2$.
For given prime ideals $P,Q\subset R$ such that $R/(P+Q)$ has finite length, let
\[
 X=\Spec R,\qquad
 Y=\Spec(R/P),\qquad Z=\Spec(R/Q).
 \]
Let $s\in X$ be the closed point.
Set $a=\dim Y$, $b=\dim Z$, and $d=\dim R$, 
thus $a+b\le d$ by Serre's dimension inequality.
By de Jong's theorem \cite[Theorems~4.1 and~6.5]{deJongAlterations}, we can choose a projective alteration $g:Z'\to Z$ where $Z'$ is integral and regular. Choose a closed embedding
$Z'\to X'=\PP^N_R$ with ideal sheaf $\III$; 
this is a regular immersion of codimension
$N+d-b$. Put
\[
 Y'=X'\times_XY,\qquad
 Z'_s=Z'\times_X\Spec k,\qquad
 E=N_{Z'/X'},\qquad E_s=E|_{Z'_s},
\]
thus $E=\Spec_{Z'}\gr_{\III}\OOO_{X'}$.
Consider the normal cone
\[
C_{Z'\cap Y'/Y'}=\Spec_E\FFF, \qquad \FFF=\gr_{\III}\OOO_{Y'},
\]
where $\FFF$ is viewed as a coherent sheaf on $E$.
Then
\begin{equation}
\label{eq:gabber-support-bound}
\dim\Supp\FFF=\dim(Y')=N+a\le N+d-b=\rank(E)
\end{equation}
where the first equality follows from \eqref{eq:cone-dimension-preliminaries} since $Y'$ acquires maximal local dimension at the closed points of $Z'_s$.
Let 
\[
\GGG_j=\m^j\FFF/\m^{j+1}\FFF,
\]
regarded as a coherent sheaf on $E_s$ for $j\ge 0$. Then $\bigoplus_j\GGG_j$ is a graded coherent algebra on $E_s$ concentrated in finitely many degrees $j$
because a power of $\m$ is contained in $P+Q$ and hence annihilates $\FFF$.
Put
\begin{equation}
\label{eq:gabber-cone-C}
C=C_{Z'\cap Y'/Y'}\times_{Z'}Z'_s=\Spec_{E_s}\GGG_0.
\end{equation}
The following diagram shows the geometric data to be used later.
\[
\xymatrix@!C@C-0.4em@M+0.2em@!C{
 C\ar[r] & E_s\ar[r]\ar[d]\ar@{}[dr]|\square & E\ar[d] \\
 & Z'_s\ar[r]\ar[d]\ar@{}[dr]|\square
 & Z'\ar[r]\ar[d]
 & X'\ar[d]\ar@{}[dr]|\square
 & Y'\ar[l]\ar[d]\ar@{}[dr]|\square
 & X'_s\ar[l]\ar[d]
 & Z'_s\ar[l] \\
 & \Spec k\ar[r] & Z\ar[r] & X & Y\ar[l] & \Spec k\ar[l]
}
\]

\subsection{Global generation}
\label{subsec:global-generation}

In the situation of \S\ref{subsec:alteration-normal-cone}
let $R_1=R/\m^2$ and $X_1=\Spec R_1$ and $X'_1=X'\times_XX_1=\PP^N_{R_1}$. Put 
\[
T=T_{\mathrm h}\oplus T_{\mathrm v}
\]
where $T_{\mathrm h}=T_sX=\Spec\Sym(\m/\m^2)$ and $T_{\mathrm v}$ is the $k$-vector space of global vector fields on $X'_s=\PP^N_k$. The derivative yields a homomorphism of locally free sheaves over $Z'_s$
\refstepcounter{equation}
\[
(\III/\III^2)|_{Z'_s}\to \Omega_{X'_1/k}|_{Z'_s}=\Omega_{X'_s/k}|_{Z'_s}\oplus(\m/\m^2)\otimes_k\OOO_{Z'_s}
\tag{\theequation}\label{eq:gabber-conormal-derivative}
\]
where $X_1'$ is a $k$-scheme by the section $k\to R/\m^2$.
Since $Z'$ and $X'$ are regular, this homomorphism is injective on the fiber at every closed point, hence locally split injective. 
Since the tangent sheaf of $\PP^N_k$ is globally generated, the dual of \eqref{eq:gabber-conormal-derivative} and evaluation of global vector fields induce a surjective homomorphism of vector bundles over $Z'_s$
\begin{equation}
\label{eq:gabber-u}
u:Z'_s\times_kT\to E_s.
\end{equation}

\subsection{Numerical consequences}

The geometric setup yields the following results. 
See \cite[Theorem~1.2]{KuranoRoberts2000} for an algebraic formulation of the positivity criterion.

\begin{proposition} 
\label{prop:gabber-numerical}
In the setting of \S\ref{subsec:alteration-normal-cone},
we have $\chi^R(R/P,R/Q)\ge0$, and
\[
 \chi^R(R/P,R/Q)>0
 \quad\Longleftrightarrow\quad
 \pr_2:u^{-1}(C)\to T\text{ is dominant},
\]
for $C$ as defined in \eqref{eq:gabber-cone-C}, and
$\chi^R(R/P,R/Q)=0$ if $a+b<d$.
\end{proposition}

\begin{proof}
Let $\delta$ be the degree of the alteration $g$. 
The projection formula gives
\[
 \chi_{X'}(\OOO_{Y'},\OOO_{Z'})
 =\sum_q(-1)^q\chi_X(\OOO_Y,R^qg_*\OOO_{Z'}),
\]
and in the Grothendieck group $G_0(Z)$ we have
\[
 \sum_q(-1)^q[R^qg_*\OOO_{Z'}]
 =\delta[\OOO_Z]+\varepsilon
\]
where $\varepsilon$ is supported in dimension less than $b$.
By induction on $b$ we may assume that 
$\chiR(R/P,R/Q')=0$ when $Q'\supset Q$ is prime with
$\dim(R/Q')<b$. It follows that
\[
 \delta\,\chiR(R/P,R/Q)
 =\chi_{X'}(\OOO_{Y'},\OOO_{Z'}).
\]
By Lemma~\ref{lem:normal-cone-euler-compare} and the finite
$\m$-adic filtration of $\FFF$ we have
\[
 \chi_{X'}(\OOO_{Y'},\OOO_{Z'})
 =\chi_E(\FFF,\OOO_{Z'})
 =\sum_{j\ge0}\chi_{E_s}(\GGG_j,\OOO_{Z'_s}).
\]
The surjective homomorphism $u$ of \eqref{eq:gabber-u} and the dimension inequality \eqref{eq:gabber-support-bound} allow us to apply
Lemma~\ref{lem:gabber-generated}. It follows that
\[
 \delta\cdot\chiR(R/P,R/Q)
 =\sum_{j\ge0}\dim_{k(T)}
 H^0\bigl((Z'_s)_{k(T)},(u^*\GGG_j)_{k(T)}\bigr).
\]
This number is always non-negative, and it is positive iff $\pr_2:u^{-1}(C)\to T$ is dominant since $C=\Spec_{E_s}\GGG_0$ and $\GGG_j$ is a $\GGG_0$-module. If $a+b<d$, the dimension inequality \eqref{eq:gabber-support-bound} is strict, so the generic fiber of $\pr_2$ is empty.
\end{proof}

\begin{theorem}[Nonnegativity and vanishing]
\label{thm:gabber-nonnegativity}
Let $M,N$ be finite modules over a regular local ring $R$ such that
$M\otimes_RN$ has finite length. Then $\chiR(M,N)\ge0$, and $\chiR(M,N)=0$ when $\dim M+\dim N<\dim R$.
\end{theorem}

\begin{proof}
By Lemma~\ref{lem:geometric-approximation-general} (approximation) we may assume that $R$ is geometric with algebraically closed residue field $k=R/\m$. In mixed characteristic $(0,p)$,
by adjoining a root of $p$ if necessary we may assume that $p\in\m^2$ so that $R$ has a first-order coefficient field. By prime filtrations we may assume that $M=R/P$ and $N=R/Q$ with $P,Q$ prime. Then Proposition~\ref{prop:gabber-numerical} yields the result.
\end{proof}

\section{The dominant cone addition criterion}
\label{sec:cone-criterion-proof}

Let again $(R,\m,k)$ be a regular local ring and $X=\Spec R$ with closed point $s\in X$. For an ideal $I\subseteq R$ we consider the cone 
\[
C_{R/I}=\Spec\gr_{\m}(R/I)
\]
as a closed subscheme of the affine space $T_{\mathrm h}=T_sX=\Spec\gr_{\m}R$. In this section we will prove Proposition \ref{thm:dominant-cone-intro} that was stated in the introduction.

\subsection{Blowup preparation and the normal line}

We use the setup of \S\ref{subsec:alteration-normal-cone} and \S\ref{subsec:global-generation} and assume in addition that $b=\dim(Z)>0$ and that the regular projective alteration $Z'\to Z$ factors through the blowup of $Z$ at the closed point, $Z'\to\Bl_sZ\to Z$.
Such alterations already occur in \cite[Proposition~3.5]{KuranoRoberts2000}.

\begin{lemma} 
\label{lem:cartier-fiber-proof}
The closed immersion $Z'_s\to Z'$ is an effective Cartier divisor,
and the closed immersion $Z'_s\to X'_s=\PP_k^N$ is regular. \end{lemma}

\begin{proof}
The ideal
$\m\OOO_Z$ becomes invertible on the blowup and stays invertible on the regular integral $Z'$, and its zero scheme is $Z'_s$.
The immersion $Z'\to X'$ is regular since
$Z'$ and $X'$ are regular, hence $Z'_s\to X'$ is regular.
Since $X'_s$ is regular, the immersion $Z'_s\to X'_s$ is regular as well, using for example \cite[Tag~09RK]{StacksProject}.
\end{proof}

We consider the following vector bundles on $Z'_s$ where $L$ is a line bundle.
\[
 L=N_{Z'_s/Z'},\qquad
 \widetilde E=N_{Z'_s/X'},\qquad E_s=N_{Z'/X'}|_{Z'_s}
\]
The regular immersions $Z'_s\to Z'\to X'$ give the
normal sequence
\[
 0\longrightarrow L\xrightarrow{\ \iota\ }\widetilde E
    \xrightarrow{\ \alpha\ }E_s\longrightarrow0.
\]
The functoriality of normal cones gives a morphism $\beta:\widetilde E\to T_{\mathrm h}$ such that the composition $\beta\circ\iota$ induces $\lambda:L\to C_{R/Q}$.
Let $\pr_{\mathrm h}:T\to T_{\mathrm h}$ be the projection.

\begin{lemma} 
\label{lem:criterion-generating-lift}
The homomorphism $u:Z'_s\times_k T\to E_s$ of \eqref{eq:gabber-u} has a lift
$\widetilde u:Z'_s\times_k T\to\widetilde E$ such that 
$\beta\circ\widetilde u={\pr_{\mathrm h}}\circ\pr_2$.
\end{lemma}

\begin{proof}
Let $\JJJ=\III+\m\OOO_{X'}$, the ideal of $Z'_s$ in $X'$. As in \S\ref{subsec:global-generation}, the derivative gives an extension of \eqref{eq:gabber-conormal-derivative} to a homomorphism 
\refstepcounter{equation}
\[
(\JJJ/\JJJ^2)|_{Z'_s}\to \Omega_{X'_1/k}|_{Z'_s}=\Omega_{X'_s/k}|_{Z'_s}\oplus(\m/\m^2)\otimes_k\OOO_{Z'_s}
\tag{\theequation}\label{eq:criterion-conormal-derivative}
\]
whose composition 
with the homomorphism $(\m/\m^2)\otimes_k\OOO_{Z'_s}\to(\JJJ/\JJJ^2)|_{Z'_s}$ induced by $X'\to X$ is the inclusion of the second factor. 
The dual of \eqref{eq:criterion-conormal-derivative} and evaluation of global vector fields give $\widetilde u$.
\end{proof}

The following diagram collects our bundle and cone morphisms so far.
\[
\xymatrix@!C@M+.2em{
 && Z'_s\times_k T \ar[d]^{\widetilde u}\ar[dr]^{u}
    \ar@/_2pc/[dd]|(.5)\hole_(.25){{\pr_{\mathrm h}}\circ\,{\pr_2}} & \\
 0\ar[r] & L\ar[r]^(.4){\iota}\ar[d]_(.45){\lambda}
    & \widetilde E\ar[r]^(.45){\alpha}\ar[d]^(.45){\beta}
    & E_s\ar[r]&0\\
 & C_{R/Q}\ar@{^(->}[r] & T_{\mathrm h} 
}
\]

\begin{lemma}
\label{lem:criterion-line}
The morphism $\lambda:L\to C_{R/Q}$ is surjective.
\end{lemma}

\begin{proof}
Let $Z''$ be the blowup of $Z$ at the closed point with exceptional fiber $Z''_s$, and $L'=N_{Z''_s/Z''}$. The alteration $g$ factors as $Z'\to Z''\to Z$, which induces morphisms 
\[
L\to L'\to C_{R/Q} \quad\text{over}\quad Z'_s\to Z''_s\to\Spec k.
\] 
Here the first map identifies $L$ with $L'\times_{Z''_s}Z'_s$. The morphism $Z'\to Z''$ is dominant and projective, hence surjective. Thus $Z'_s\to Z''_s$ and hence $L\to L'$ are surjective. The map $L'\to C_{R/Q}$ is surjective by standard properties of the blowup; it restricts to an isomorphism of the complements of the zero sections.
\end{proof}

\begin{lemma} 
\label{lem:criterion-fiber-cone}
The homomorphisms $\alpha$ and $\beta$ induce the following commutative diagram where {\footnotesize $\square$} is Cartesian and $C$ has been defined in \eqref{eq:gabber-cone-C}.
\[
\xymatrix@M+0.2em{
 C_{R/P}\ar@{^(->}[d] & 
 C_{Z'_s/Y'}\ar[l]\ar[r]\ar@{^(->}[d]\ar@{}[dl]|\square & 
 C\ar@{^(->}[d]\\ 
 T_{\mathrm h} & \widetilde E\ar[l]_-{\beta}\ar[r]^-{\alpha} & E_s .
}
\]
In other words, $\alpha$ induces a morphism $\beta^{-1}(C_{R/P})\to C$.
\end{lemma}

\begin{proof}
Lemma~\ref{lem:cone-comparison-preliminaries} gives the Cartesian square on the left. Functoriality of the normal cone induces $C_{Z'_s/Y'}\to C_{Z'\cap Y'/Y'}$ and hence the right hand square.
\end{proof}

\subsection{The fundamental diagram}

The idea is to use the freedom in choosing a lift of a vector
in $E_s$ to $\widetilde E$ to prescribe its image in $T_{\mathrm h}$.

\begin{proposition}
\label{prop:criterion-incidence-direct}
There is a commutative diagram
\[
\xymatrix@M+0.2em{
 L\times_k C_{R/P}\times_k T_{\mathrm v}
   \ar[r]^-{\widetilde\alpha}\ar[d]_{\lambda\times\id\times\id} &
 u^{-1}(C)\ar[d]^{\pr_2}\\
 C_{R/Q}\times_k C_{R/P}\times_k T_{\mathrm v}
   \ar[r]^-{+} & T .
}
\]
Here $+$ is the addition in $T$, and $\widetilde\alpha$ is the
restriction of the morphism
\[
 \widetilde\alpha:L\times_k T\longrightarrow Z'_s\times_k T,
 \qquad (\ell,t)\longmapsto(\pi(\ell),t+\lambda(\ell)),
\]
where $\pi:L\to Z'_s$ is the projection and $\lambda(\ell)$ is
viewed as an element of $T_{\mathrm h}\subseteq T$.
\end{proposition}

\begin{proof}
We have to show that $u\circ\widetilde\alpha$ maps $C'=L\times_k C_{R/P}\times_k T_{\mathrm v}$ into $C$.
If this holds, the diagram commutes by the definition of $\widetilde\alpha$.
One verifies directly that the following diagram commutes, where $v=\widetilde u\circ\widetilde\alpha-\iota\circ\pr_1$.
\[
\xymatrix@C+2em@M+0.2em{
 T\ar[d]_{\pr_{\mathrm h}}
 & L\times_k T
     \ar[l]_-{\pr_2}\ar[r]^-{\widetilde\alpha}\ar[d]^v
 & Z'_s\times_k T\ar[d]^u\\
 T_{\mathrm h}
 & \widetilde E\ar[l]_-{\beta}\ar[r]^-{\alpha}
 & E_s
}
\]
Indeed, the relations $\alpha\circ\widetilde u=u$ and $\beta\circ\widetilde u={\pr_{\mathrm h}}\circ\pr_2$ of Lemma \ref{lem:criterion-generating-lift} together with $\alpha\circ\iota=0$ and $\beta\circ\iota=\lambda$
give
\[
\begin{aligned}
 \alpha\circ v
 &=\alpha\circ(\widetilde u\circ\widetilde\alpha-\iota\circ\pr_1)
 =u\circ\widetilde\alpha,\\
 \beta\circ v
 &={\pr_{\mathrm h}}\circ{\pr_2}\circ\widetilde\alpha
   -\lambda\circ\pr_1
 ={\pr_{\mathrm h}}\circ\pr_2.
\end{aligned}
\]
Clearly the morphism $\beta\circ v=\pr_{\mathrm h}\circ\pr_2$ maps $C'$ into $C_{R/P}$. Hence $v$ maps $C'$ into $C_{Z'_s/Y'}$, and $u\circ\widetilde\alpha=\alpha\circ v$ maps $C'$ into $C$ as required.
\end{proof}

\begin{proposition}
\label{prop:cone-criterion-proof}
Assume that $R$ is geometric with first-order coefficient field and algebraically closed residue field. Let
$I,J\subseteq R$ be ideals such that $R/(I+J)$ has finite length.  
If the addition map
\[
 C_{R/I}\times_k C_{R/J}\longrightarrow T_R,\qquad
 (x,y)\longmapsto x+y
\]
is dominant, then $\chiR(R/I,R/J)>0$.
\end{proposition}

\begin{proof}
Dominance implies that $\dim R/I+\dim R/J=\dim C_{R/I}+\dim C_{R/J}\ge\dim R$, so equality holds by Serre's dimension inequality.
We may assume that $R/J$ has positive dimension.
By nonnegativity (Proposition~\ref{prop:gabber-numerical}) and Lemma~\ref{lem:cone-minimal-primes-preliminaries} we can assume that $I=P$ and $J=Q$ are prime.
Then we have the diagram of Proposition~\ref{prop:criterion-incidence-direct} where $\lambda$ is surjective by Lemma~\ref{lem:criterion-line} and the addition map is dominant. It follows that $\pr_2$ is dominant, and the assertion follows by Proposition~\ref{prop:gabber-numerical}.
\end{proof}

\section{Transversality via symmetric covers}
\label{sec:symmetric-covers-direct}

For a field $k$ and $n\ge2$, let $V_n$ be the kernel of the addition map $k^n\to k$.
The symmetric group $\mathfrak S_n$ acts on $V_n$ by the permutation
action on $k^n$.

\begin{lemma} 
\label{lem:symmetric-transversality-direct}
Assume that $n\ge 2$ is invertible in $k$. Then for subspaces
$U,W\subseteq V_n$ with $\dim U+\dim W=\dim V_n$ there is an element
$g\in\mathfrak S_n$ with $gU\cap W=0$.
\end{lemma}

\begin{proof}
For a non-zero $v\in V_n$, the coordinates of $v$ are not all equal, and
\[
 \Span_k\{(\tau-1)v\mid \tau\text{ a transposition}\}=V_n
\]
since $(\tau-1)v=(v_j-v_i)(e_i-e_j)$ for $\tau=(ij)$,
and the graph joining $i$ to $j$ when $v_i\ne v_j$ is connected.
Choose a translate of $U$ minimizing $\dim(U\cap W)$. 
If $v\in U\cap W$ is non-zero, choose a transposition $\tau$ with $(\tau-1)v$ outside $U+W$. For $x\in U$ with $\tau x\in W$, the difference $\tau x-x$ belongs to $U+W$ and
to the line spanned by $(\tau-1)v$, hence is zero. Thus $\tau U\cap W$ is the $\tau$-fixed subspace in
$U\cap W$. It excludes $v$, contradicting minimality.
\end{proof}

\begin{lemma} 
\label{lem:invariant-cones-direct}
Let $n\ge2$ be invertible in $k$ and
view $V_n$ as an affine space over $k$.
For closed cones $C_1,C_2\subseteq V_n$ where
$C_1$ is $\mathfrak S_n$-invariant and
$\dim C_1+\dim C_2=\dim V_n$,
 the addition map
$C_1\times_kC_2\to V_n$ is dominant.
\end{lemma}

\begin{proof}
We may assume that $k$ is algebraically closed and that the $C_i$ are reduced. Let $x_i\in C_i$ be smooth $k$-points on irreducible components of maximal dimension.
The tangent spaces $U_i=T_{C_i,x_i}$ are subspaces of $V_n$ with complementary dimensions.
Lemma~\ref{lem:symmetric-transversality-direct} gives a $g\in\mathfrak S_n$ with
$gU_1\cap U_2=0$. 
The addition map is \'etale at $(gx_1,x_2)$ and hence dominant because its differential is an isomorphism.
\end{proof}

\begin{proposition} 
\label{lem:symmetric-cover-direct}
Let $R$ be a regular local ring of positive dimension $d$, and put $n=d+1$. There is a finite free regular local
$R$-algebra $S$ with the same residue field and an action of $\mathfrak S_n$ by $R$-algebra automorphisms such that the tangent representation $T_{S}$ is isomorphic to $V_n$ and $S$ has a first-order coefficient field.
\end{proposition}

\begin{proof}
Let $s_i\in\ZZ[t_1,\ldots,t_n]$ be the $i$-th elementary
symmetric polynomial. The $s_i$ are algebraically independent, and
$
 \ZZ[s_1,\ldots,s_n]\to\ZZ[t_1,\ldots,t_n]
$
is finite free of rank $n!$ by
\cite[Ch.~IV, \S6, no.~1, Theorem~1]{BourbakiAlgebra}.
Choose regular parameters $x_1,\ldots,x_d$ of $R$ and form
\[
 S=R\otimes_{\ZZ[s_1,\ldots,s_n]}\ZZ[t_1,\ldots,t_n]
\]
using $s_1\mapsto0$ and $s_{i+1}\mapsto x_i$.
Then $S$ has the desired properties.
Indeed, the image of $t_i$ is nilpotent in 
$S/\m_R S$, hence $S$ is local with
residue field $R/\m_R$ and maximal ideal $\m_{S}=(t_1,\ldots,t_n)$.
The identities $x_i=s_{i+1}(t)$ give $\m_R S\subseteq\m_{S}^2$, 
and $s_1(t)=0$ implies that $\m_S$ is generated by $d$ elements. 
Moreover $\dim S=d$, hence $S$ is regular. The group $\mathfrak S_n$ acts by permutations of the $t_i$, and the vector space $\m_S/\m_S^2$ is spanned by the $t_i$ with sole relation $\sum t_i=0$. All assertions follow.
\end{proof}

The following result was originally proved in \cite{OpenAI2026} using different methods.

\begin{theorem} 
\label{thm:positivity-direct}
Let $M,N$ be finite modules over a regular local ring $R$ such that $M\otimes_R N$ has finite
length and $\dim M+\dim N=\dim R$. Then $\chi^R(M,N)>0$.
\end{theorem}

\begin{proof}
By Lemma~\ref{lem:geometric-approximation-general} we may
assume that $R$ is geometric with algebraically closed residue
field. Put $d=\dim R$.
By nonnegativity (Theorem~\ref{thm:gabber-nonnegativity})
it suffices to consider $M=R/I$ and $N=R/J$
with $\dim R/I+\dim R/J=d$.

By adjoining a variable $t$, localizing at $(\m_R,t)$, and letting $t$ act trivially on $N$ we may increase $d$ and hence assume that $n=d+1$ is invertible in $R$ and $n\ge 2$.
By Proposition~\ref{lem:symmetric-cover-direct} and flat base change
we may assume
that the symmetric group $\mathfrak S_n$ acts on $R$ such that
$I$ and $J$ are invariant, $T_R\cong V_n$ as a representation of $\mathfrak S_n$, and $R$ has a first-order coefficient field.
Then the addition map $C_{R/I}\times C_{R/J}\to T_R$ is dominant
by Lemma~\ref{lem:invariant-cones-direct}, and
Proposition~\ref{prop:cone-criterion-proof} gives
$\chi^R(R/I,R/J)>0$.
\end{proof}

\begin{example}
Let us trace the effect of the symmetric cover in an example. 
Let $R=k[x,y]_{(x,y)}$, $I=(y)$ and $J=(y-x^2)$.
The tangent cones of $R/I$ and $R/J$ are both the line $y=0$.
For the symmetric cover
$S=k[t_1,t_2,t_3]_{(t_1,t_2,t_3)}/(s_1)$, with $x=s_2$ and $y=s_3$,
the ideals $IS$ and $JS$ have the same initial equation
$t_1t_2t_3=0$. Their tangent cones are each a union of three distinct
lines in $T_S=V_3$, and their sum is $T_S$.
For two distinct lines, one from
each tangent cone, the sum is the tangent plane.
\end{example}

\end{document}